\documentclass[12pt]{amsart}

\newtheorem{theorem}{Theorem}[section]
\newtheorem{lemma}[theorem]{Lemma}
\newtheorem{corollary}[theorem]{Corollary}

\newtheorem{conjecture}[theorem]{Conjecture}

\theoremstyle{plain}
\newtheorem{definition}[theorem]{Definition}

\newtheorem{remark}[theorem]{Remark}

\theoremstyle{definition}

\theoremstyle{remark}

\numberwithin{equation}{section}

\usepackage{CJK}
\usepackage{fancyhdr}
\usepackage{amscd}
\usepackage{amsmath}
\usepackage{amsthm}
\usepackage{amssymb}
\begin{document}

\title[]{Infinitesimal Deformation of Deligne cycle class map}

\author{Sen Yang}
\address{Shing-Tung Yau Center of Southeast University \\ 
Southeast University \\
Nanjing, China\\
}

\address{School of Mathematics \\  Southeast University \\
Nanjing, China\\
}
\address{Yau Mathematical Sciences Center \\
Tsinghua University \\
Beijing, China\\
}
\email{101012424@seu.edu.cn; syang@math.tsinghua.edu.cn}

\date{}

\maketitle

\begin{abstract}
 In this note, we study the infinitesimal forms of Deligne cycle class maps. As an application, we prove that the infinitesimal form of a conjecture by Beilinson \cite{Beilinson} is true.
\end{abstract}

\tableofcontents

\section{Introduction}
\label{Introduction}
In \cite{Beilinson}, Beilinson made the following conjecture:
\begin{conjecture} [Conjecture 2.4.2.1 \cite{Beilinson}] \label{conjecture: main}
Let X be a smooth projective variety defined over a number field
k, then for each positive integer $p$, the rational Chow group $CH^{p}(X)_{\mathbb{Q}}$ injects into Deligne cohomology of  $X_{\mathbb{C}}$, where $X_{\mathbb{C}}:= X \times_{k} \mathbb{C}$. Concretely,
if a class in $CH^{p}(X)_{\mathbb{Q}}$ vanishes in Deligne cohomology  $H_{\mathcal{D}}^{2p}(X_{\mathbb{C}}, \mathbb{Z}(p))_{\mathbb{Q}}$ under the composition
\[
CH^{p}(X)_{\mathbb{Q}} \rightarrow CH^{p}(X_{\mathbb{C}})_{\mathbb{Q}} \xrightarrow{r}  H_{\mathcal{D}}^{2p}(X_{\mathbb{C}}, \mathbb{Z}(p))_{\mathbb{Q}},
\]
where the right arrow $r$ is the cycle class map for Deligne cohomology, then it is 0. 

\end{conjecture}
This conjecture is very difficult to approach, and up to now {\it there is not a single example with dimension $X \geq 2$ with large Chow ring $CH(X \times_{k} \mathbb{C})$ for which this conjecture has been verified}.\footnote{See the first paragraph of  page 1 \cite{EH}.} Esnault and Harris \cite{EH} suggest a modest conjecture (see Theorem 0.1 \cite{EH}) which follows from Conjecture 1.1 and has been proved in a particular case (see Theorem 0.2 \cite{EH}) by using $l$-adic cohomology.

The main result of this note is to study the infinitesimal form of the Deligne cycles class map, see Lemma \ref{lemma: inf'l R-first} and Theorem \ref{theorem: inf'l form2}. As an application, we prove that the infinitesimal form of the Conjecture \ref{conjecture: main}
is true, see Theorem \ref{theorem: main theorem}.

In a companion paper \cite{Yang}, as a further application, we show that the infinitesimal form of the following conjecture (due to Griffiths-Harris) is true, 
\begin{conjecture}[\cite{GH}] 
Let $X \subset \mathbb{P}_{\mathbb{C}}^{4}$ be a general hypersurface of degree $d \geq 6$, we use 
\[
\psi: CH_{hom}^{2}(X) \to J^{2}(X)
\]
to denote the Abel-Jacobi map from algebraic 1-cycles on $X$ homologically equivalent to zero to the intermediate Jacobian $J^{2}(X)$, $\psi$ is zero.
\end{conjecture}

\section{Main results}
Let $X$ be a smooth projective variety defined over $\mathbb{C}$ and $\mathbb{C}[\varepsilon]/(\varepsilon^{2})$ be the ring of dual numbers, we use $X[\varepsilon]:= X \times \mathrm{Spec}(\mathbb{C}[\varepsilon]/(\varepsilon^{2}))$ to denote the first order infinitesimal deformation of $X$. 
The classical definition of Chow groups can not recognize nilpotents, so to overcome this deficiency, for each positive integer $p$, one uses the following Soul\'e's variant of Bloch-Quillen identification to study the infinitesimal deformations of  Chow groups,
\begin{equation}
 CH^{p}(X)= H^{p}(X, K^{M}_{p}(O_{X})) \ \mathrm{ modulo \ torsion},
\end{equation}
where $K^{M}_{p}(O_{X})$ is the 
Milnor K-theory sheaf associated to the presheaf
\[
  U \to K^{M}_{p}(O_{X}(U)).
\]
Using the identification (2.1),  one considers $H^{p}(X, K^{M}_{p}(O_{X[\varepsilon]}))$ as the first order infinitesimal deformation of  $CH^{p}(X)$
 and defines, 
 \begin{definition} \label{definition: tangentChow}
 Let $X$ be a smooth projective variety defined over $\mathbb{C}$, for each positive integer $p$, the formal tangent space to $CH^{p}(X)$, denoted $T_{f}CH^{p}(X)$,  is defined to be the kernel of the natural map
\[
H^{p}(X, K^{M}_{p}(O_{X[\varepsilon]})) \xrightarrow{\varepsilon =0} H^{p}(X, K^{M}_{p}(O_{X})).
\]
\end{definition}
It is known that the formal tangent space $T_{f}CH^{p}(X)$ can be identified with $H^{p}(X, \Omega_{X/ \mathbb{Q}}^{p-1})$, where $\Omega_{X/ \mathbb{Q}}^{p-1}$ is the absolute differential.

Deligne cohomology $H_{\mathcal{D}}^{2p}(X, \mathbb{Z}(p))$ is defined to be the hypercohomology of the Deligne complex $\mathbb{Z}(p)_{\mathcal{D}}$ (in analytic topology) :
\[
\mathbb{Z}(p)_{\mathcal{D}}: \ 0 \to \mathbb{Z}(p) \to O_{X} \to \cdots \to \Omega_{X/ \mathbb{C}}^{p-1} \to 0,
\]
where $\mathbb{Z}(p)=(2 \pi i)^{p}\mathbb{Z}$ is in degree 0.
The infinitesimal deformation of this complex, denoted $\mathbb{Z}(p)_{\mathcal{D}}[\varepsilon]$, has the form, 
\[
\mathbb{Z}(p)_{\mathcal{D}}[\varepsilon]: \ 0 \to \mathbb{Z}(p) \to O_{X[\varepsilon]} \to \cdots \to \Omega_{X[\varepsilon]/ \mathbb{C}[\varepsilon]}^{p-1} \to 0,
\]
where $\mathbb{Z}(p)$ is still equal to $(2 \pi i)^{p}\mathbb{Z}$.
\begin{definition}
Let $X$ be a smooth projective variety defined over $\mathbb{C}$, for each positive integer $p$, the tangent complex to the Deligne complex $\mathbb{Z}(p)_{\mathcal{D}}$, 
denoted $\overline{\mathbb{Z}}(p)_{\mathcal{D}}$, is defined to be the kernel of  the natural map
\[
\mathbb{Z}(p)_{\mathcal{D}}[\varepsilon]  \xrightarrow{\varepsilon =0} \mathbb{Z}(p)_{\mathcal{D}}.
\]
\end{definition}

Since the map $X \to X[\varepsilon]$ has a retraction $X[\varepsilon] \to X$, $\Omega_{X[\varepsilon]/ \mathbb{C}[\varepsilon]}^{i} = \Omega_{X/ \mathbb{C}}^{i} \oplus \varepsilon \Omega_{X/ \mathbb{C}}^{i}$, where $i=0, \cdots, p-1$. The tangent complex to the Deligne complex is a direct summand of the thickened Deligne complex 
\[
\mathbb{Z}(p)_{\mathcal{D}}[\varepsilon]  = \mathbb{Z}(p)_{\mathcal{D}} \oplus \overline{\mathbb{Z}}(p)_{\mathcal{D}}.
\]
One easily sees that 
\begin{lemma}
Let $X$ be a smooth projective variety defined over $\mathbb{C}$, for each positive integer $p$, the tangent complex (to the Deligne complex $\mathbb{Z}(p)_{\mathcal{D}}$) $\overline{\mathbb{Z}}(p)_{\mathcal{D}}$ is of the form
\[
0 \to O_{X} \to \cdots \to \Omega_{X/ \mathbb{C}}^{p-1} \to 0,
\]
where $O_{X}$ is in degree 1 and $\Omega_{X/ \mathbb{C}}^{p-1}$ is in degree $p$.

\end{lemma}

We consider the hypercohomology $\mathbb{H}^{2p}(X, \mathbb{Z}(p)_{\mathcal{D}}[\varepsilon])$ of the complex $\mathbb{Z}(p)_{\mathcal{D}}[\varepsilon]$ as the infinitesimal deformation of the Deligne cohomology $H_{\mathcal{D}}^{2p}(X, \mathbb{Z}(p))$ and define,
\begin{definition}
Let $X$ be a smooth projective variety defined over $\mathbb{C}$, for each positive integer $p$, the formal tangent space to the Deligne cohomology
$H_{\mathcal{D}}^{2p}(X, \mathbb{Z}(p))$, denoted $T_{f}H_{\mathcal{D}}^{2p}(X, \mathbb{Z}(p))$, is defined to be the kernel of  the natural map 
\[
\mathbb{H}^{2p}(X, \mathbb{Z}(p)_{\mathcal{D}}[\varepsilon]) \xrightarrow{\varepsilon =0} \mathbb{H}^{2p}(X, \mathbb{Z}(p)_{\mathcal{D}}).
\]
\end{definition}

By the definition, the formal tangent space to the Deligne cohomology is the hypercohomology of the tangent complex $\overline{\mathbb{Z}}(p)_{\mathcal{D}}$:
\[
T_{f}H_{\mathcal{D}}^{2p}(X, \mathbb{Z}(p)) = \mathbb{H}^{2p}(X, \overline{\mathbb{Z}}(p)_{\mathcal{D}}) = \mathbb{H}^{2p-1}(X, \overline{\mathbb{Z}}(p)_{\mathcal{D}}[-1]).
\]
The following isomorphism is a standard fact in complex geometry,
\begin{lemma} [cf. Proposition on page 17 of \cite{Green-note}]With the notations above, one has the isomorphism
\[
\mathbb{H}^{2p}(X, \overline{\mathbb{Z}}(p)_{\mathcal{D}}) \cong 
H^{2p-1}(O_{X})\oplus H^{2p-2}(\Omega_{X/ \mathbb{C}}^{1}) \oplus \cdots \oplus H^{p}(\Omega_{X/ \mathbb{C}}^{p-1}).
\]
\end{lemma}

Let $\overline{K}^{M}_{p}(O_{X})$ denote the kernel of  the natural map 
\[
 K^{M}_{p}(O_{X[\varepsilon]}) \xrightarrow{\varepsilon =0}  K^{M}_{p}(O_{X}).
\]
Since the map $X \to X[\varepsilon]$ has a retraction $X[\varepsilon] \to X$, $K^{M}_{p}(O_{X[\varepsilon]}) = K^{M}_{p}(O_{X}) \oplus \overline{K}^{M}_{p}(O_{X})$. By Definition \ref{definition: tangentChow}, the tangent space $T_{f}CH^{p}(X)$ is $H^{p}(X, \overline{K}^{M}_{p}(O_{X}))$. Next, we would like to construct a map between tangent spaces
\[
H^{p}(X, \overline{K}^{M}_{p}(O_{X}))  \to \mathbb{H}^{2p}(X, \overline{\mathbb{Z}}(p)_{\mathcal{D}}),
\]
which is the infinitesimal form of the Deligne cycle class map 
$CH^{p}(X) \to H_{\mathcal{D}}^{2p}(X, \mathbb{Z}(p))$.

An element of  $\overline{K}^{M}_{p}(O_{X})$ is of the form 
$\prod_{i} \{f^{i}_{1}+ \varepsilon g^{i}_{1}, \cdots,  f^{i}_{p}+ \varepsilon g^{i}_{p}  \} \in K^{M}_{p}(O_{X[\varepsilon]})$ such that $\prod_{i} \{f^{i}_{1}, \cdots,  f^{i}_{p} \} =1 \in K^{M}_{p}(O_{X})$. We are reduced to looking at
$\{f_{1}+ \varepsilon g_{1}, \cdots,  f_{p}+ \varepsilon g_{p}  \} \in K^{M}_{p}(O_{X[\varepsilon]})$ such that $\{f_{1}, \cdots, f_{p}   \}=1 \in K^{M}_{p}(O_{X})$. For simplicity, we assume that $g_{2}= \cdots=g_{p}=0$ and have
\begin{align*}
\{f_{1}+ \varepsilon g_{1}, f_{2}, \cdots, f_{p}  \} & = \{f_{1}, f_{2}, \cdots, f_{p}  \} \{1+ \varepsilon \dfrac{g_{1}}{f_{1}}, f_{2}, \cdots, f_{p}  \}  \\
& =  \{1+ \varepsilon \dfrac{g_{1}}{f_{1}}, f_{2}, \cdots, f_{p}  \}.
\end{align*}

Applying $\wedge^{p} dlog$ to $\{1+ \varepsilon \dfrac{g_{1}}{f_{1}}, f_{2}, \cdots, f_{p}  \}$, where $d= d_{\mathbb{C}[\varepsilon]}$, one obtains 
\begin{align*}
\dfrac{d_{\mathbb{C}[\varepsilon]}(1+ \varepsilon\dfrac{g_{1}}{f_{1}})}{1+ \varepsilon\dfrac{g_{1}}{f_{1}}} \wedge \dfrac{d_{\mathbb{C}}f_{2}}{f_{2}} \wedge \cdots \wedge  \dfrac{d_{\mathbb{C}}f_{p}}{f_{p}} & = (1- \varepsilon\dfrac{g_{1}}{f_{1}})d_{\mathbb{C}[\varepsilon]}(1+ \varepsilon\dfrac{g_{1}}{f_{1}}) \wedge \dfrac{d_{\mathbb{C}}f_{2}}{f_{2}} \wedge \cdots \wedge  \dfrac{d_{\mathbb{C}}f_{p}}{f_{p}} \\
& = \varepsilon   d_{\mathbb{C}}(\dfrac{g_{1}}{f_{1}}) \wedge \dfrac{d_{\mathbb{C}}f_{2}}{f_{2}} \wedge \cdots \wedge  \dfrac{d_{\mathbb{C}}f_{p}}{f_{p}} \\
&  = \varepsilon \dfrac{f_{1}d_{\mathbb{C}} g_{1}- g_{1}d_{\mathbb{C}}f_{1}}{f^{2}_{1}} \wedge \dfrac{d_{\mathbb{C}}f_{2}}{f_{2}} \wedge \cdots \wedge  \dfrac{d_{\mathbb{C}}f_{p}}{f_{p}}. 
\end{align*}
This gives a map from $\overline{K}^{M}_{p}(O_{X})$  to $  \Omega_{X/ \mathbb{C}}^{p}$,
\begin{definition}
One defines a map $ \tilde{\wedge}^{p} dlog: \overline{K}^{M}_{p}(O_{X}) \to \Omega_{X/ \mathbb{C}}^{p}$ by
\begin{equation}
 \{1+ \varepsilon \dfrac{g_{1}}{f_{1}}, f_{2}, \cdots, f_{p}  \} \to \dfrac{f_{1}d_{\mathbb{C}}g_{1}- g_{1}d_{\mathbb{C}}f_{1}}{f^{2}_{1}} \wedge \dfrac{d_{\mathbb{C}}f_{2}}{f_{2}} \wedge \cdots \wedge  \dfrac{d_{\mathbb{C}}f_{p}}{f_{p}}.
\end{equation}
\end{definition}

On the other hand, let 
\begin{equation}
\omega= \dfrac{g_{1}}{f_{1}} \dfrac{d_{\mathbb{C}}f_{2}}{f_{2}} \wedge \cdots \wedge  \dfrac{d_{\mathbb{C}}f_{p}}{f_{p}} \in  \Omega_{X/ \mathbb{C}}^{p-1},
\end{equation}
 one checks that 
\begin{align*}
d_{\mathbb{C}}(\omega) & = d_{\mathbb{C}}(\dfrac{g_{1}}{f_{1}f_{2} \cdots f_{p}}) \wedge d_{\mathbb{C}}f_{2}\wedge \cdots \wedge  d_{\mathbb{C}}f_{p} \\
& = \dfrac{(f_{1}f_{2} \cdots f_{p})d_{\mathbb{C}}g_{1}-g_{1}d_{\mathbb{C}} (f_{1}f_{2} \cdots f_{p})} {(f_{1}f_{2} \cdots f_{p})^{2}}\wedge d_{\mathbb{C}}f_{2} \wedge \cdots \wedge  d_{\mathbb{C}}f_{p} \\
& = (\dfrac{d_{\mathbb{C}}g_{1}} {f_{1}f_{2} \cdots f_{p}} - \dfrac{g_{1}f_{2} \cdots f_{p} d_{\mathbb{C}}f_{1}}{(f_{1}f_{2} \cdots f_{p})^{2}})\wedge d_{\mathbb{C}}f_{2} \wedge \cdots \wedge  d_{\mathbb{C}}f_{p}.
\end{align*}

Comparing $d_{\mathbb{C}} (\omega)$ with (2.2), one sees that
\begin{lemma}
With the notations above, one has
\[
 \tilde{\wedge}^{p} dlog (\{1+ \varepsilon \dfrac{g_{1}}{f_{1}}, f_{2}, \cdots, f_{p}  \}) = d_{\mathbb{C}} (\omega),
\]
where $\omega= \dfrac{g_{1}}{f_{1}} \dfrac{d_{\mathbb{C}}f_{2}}{f_{2}} \wedge \cdots \wedge  \dfrac{d_{\mathbb{C}}f_{p}}{f_{p}}$.
\end{lemma}

The following commutative diagram, which gives a quasi-isomorphism from the upper complex to the bottom one, is the tangent to the commutative diagram in Section 2.7 \cite{EV} (page 56),
\[
\begin{CD}
      O_{X} @>>>  \cdots  @>>> \Omega_{X/ \mathbb{C}}^{p-2}
      @>>> \Omega_{X/ \mathbb{C}}^{p-1}  @>>>  0 \\
      @V\alpha_{1}VV  @VVV  @V\alpha_{p-1}VV  @V \alpha_{p}VV  @VVV \\
     O_{X} @>-\delta_{1}>>  \cdots  @>-\delta_{p-2}>> \Omega_{X/ \mathbb{C}}^{p-2}
      @>-\delta_{p-1}>> \Omega_{X/ \mathbb{C}}^{p} \oplus \Omega_{X/ \mathbb{C}}^{p-1}  @>-\delta_{p}>> \Omega_{X/ \mathbb{C}}^{p+1} \oplus \Omega_{X/ \mathbb{C}}^{p} \cdots,
  \end{CD}
\]
where $\alpha_{i}(x)=(-1)^{i-1}(x)$ for $1\leq i \leq p-1$ and $\alpha_{p}(x)=(-1)^{p-1}(d_{\mathbb{C}}x, x)$; $\delta_{i}(x)=dx$ for $1\leq i \leq p-2$, 
$\delta_{p-1}(x)=(0, dx)$ and $\delta_{p}(x, y)=(-dx, -x+dy)$. 
For $(d_{\mathbb{C}}\omega, \omega) \in \Omega_{X/ \mathbb{C}}^{p} \oplus \Omega_{X/ \mathbb{C}}^{p-1}$, there exists the unique $(-1)^{p-1} \omega \in \Omega_{X/ \mathbb{C}}^{p-1}$ such that $\alpha_{p}((-1)^{p-1} \omega ) = (d_{\mathbb{C}}\omega, \omega)$.

\begin{definition}
One defines a map $\beta: \overline{K}^{M}_{p}(O_{X}) \longrightarrow \Omega_{X/ \mathbb{C}}^{p-1}$ by 
\begin{equation}
\{1+ \varepsilon \dfrac{g_{1}}{f_{1}}, f_{2}, \cdots, f_{p} \}  \to  (-1)^{p-1}\omega,
\end{equation}
Where $\omega=\dfrac{g_{1}}{f_{1}} \dfrac{d_{\mathbb{C}}f_{2}}{f_{2}} \wedge \cdots \wedge  \dfrac{d_{\mathbb{C}}f_{p}}{f_{p}}$.
\end{definition}

We can see the map $\beta$ (2.4) in an alternative way \footnote{We thank Spencer Bloch and Jerome Hoffman for comments}. Firstly, applying $\wedge^{p} dlog$ to $\{1+ \varepsilon \dfrac{g_{1}}{f_{1}}, f_{2}, \cdots, f_{p}  \}$, where $d= d_{\mathbb{C}}$, one obtains 
\begin{align*}
\dfrac{d_{\mathbb{C}}(1+ \varepsilon\dfrac{g_{1}}{f_{1}})}{1+ \varepsilon\dfrac{g_{1}}{f_{1}}} \wedge \dfrac{d_{\mathbb{C}}f_{2}}{f_{2}} \wedge \cdots \wedge  \dfrac{d_{\mathbb{C}}f_{p}}{f_{p}} & = (1- \varepsilon\dfrac{g_{1}}{f_{1}})d_{\mathbb{C}}(1+ \varepsilon\dfrac{g_{1}}{f_{1}}) \wedge \dfrac{d_{\mathbb{C}}f_{2}}{f_{2}} \wedge \cdots \wedge  \dfrac{d_{\mathbb{C}}f_{p}}{f_{p}} \\
& = (\varepsilon d_{\mathbb{C}}(\dfrac{g_{1}}{f_{1}})+ \dfrac{g_{1}}{f_{1}} d_{\mathbb{C}}(\varepsilon) ) \wedge \dfrac{d_{\mathbb{C}}f_{2}}{f_{2}} \wedge \cdots \wedge  \dfrac{d_{\mathbb{C}}f_{p}}{f_{p}};
\end{align*}
Secondly, applying the truncation map $\rfloor  \dfrac{\partial}{\partial \varepsilon}\mid_{\varepsilon =0}: \Omega_{X[\varepsilon]/\mathbb{C}}^{p} \to \Omega_{X/\mathbb{C}}^{p-1}$ to $(\varepsilon d_{\mathbb{C}}(\dfrac{g_{1}}{f_{1}})+ \dfrac{g_{1}}{f_{1}} d_{\mathbb{C}}(\varepsilon) ) \wedge \dfrac{d_{\mathbb{C}}f_{2}}{f_{2}} \wedge \cdots \wedge  \dfrac{d_{\mathbb{C}}f_{p}}{f_{p}}$, one obtains $(-1)^{p-1} \dfrac{g_{1}}{f_{1}} \dfrac{d_{\mathbb{C}}f_{2}}{f_{2}} \wedge \cdots \wedge  \dfrac{d_{\mathbb{C}}f_{p}}{f_{p}}$. This shows that the composition
\[
\overline{K}^{M}_{p}(O_{X}) \xrightarrow{\wedge^{p} dlog} \Omega_{X[\varepsilon]/\mathbb{C}}^{p} \xrightarrow{\rfloor  \dfrac{\partial}{\partial \varepsilon}\mid_{\varepsilon =0}} \Omega_{X/\mathbb{C}}^{p-1},
\]
agrees with the map $\beta$ (2.4).

The map $\beta$ (2.4) induces a map from the complex
\[
0 \to 0 \to \cdots \to 0 \to \overline{K}^{M}_{p}(O_{X}),
\]
where $\overline{K}^{M}_{p}(O_{X})$ is in degree p, to the complex $\overline{\mathbb{Z}}(p)_{\mathcal{D}}$,
\begin{equation}
\begin{CD}
      0  @>>>  0 @>>> \cdots  @>>> 0 @>>> \overline{K}^{M}_{p}(O_{X})   \\
      @V0VV  @V0VV  @V0VV  @V0VV @V \beta V(2.4)V   \\
     0 @>>> O_{X} @>>>  \cdots  @>>> \Omega_{X/ \mathbb{C}}^{p-2}   @>>> \Omega_{X/ \mathbb{C}}^{p-1}.
  \end{CD}
\end{equation}

This induces a map between (hyper)cohomology groups
 \begin{equation} 
\lambda: H^{p}(X, \overline{K}^{M}_{p}(O_{X}))  \to \mathbb{H}^{2p}(X, \overline{\mathbb{Z}}(p)_{\mathcal{D}}),
\end{equation}
which is the infinitesimal form of the Deligne cycle class map
\[
r: CH^{p}(X) \to H_{\mathcal{D}}^{2p}(X, \mathbb{Z}(p)).
\]

By Lemma 2.4, $\mathbb{H}^{2p}(X, \overline{\mathbb{Z}}(p)_{\mathcal{D}}) \cong H^{2p-1}(O_{X})\oplus \cdots \oplus H^{p}(\Omega_{X/ \mathbb{C}}^{p-1})$ and by the diagram (2.5), we note the image of 
the map (2.6) lies in $H^{p}(\Omega_{X/ \mathbb{C}}^{p-1})$. So the map (2.6) is indeed the composition:
\[
H^{p}(X, \overline{K}^{M}_{p}(O_{X})) \to H^{p}(\Omega_{X/ \mathbb{C}}^{p-1}) \hookrightarrow \mathbb{H}^{2p}(X, \overline{\mathbb{Z}}(p)_{\mathcal{D}}) .
\]

In summary,
\begin{lemma} \label{lemma: inf'l R-first}
Let $X$ be a smooth projective variety defined over $\mathbb{C}$, for each positive integer $p$, the infinitesimal form of the Deligne cycle class map
\[
r: CH^{p}(X) \to H_{\mathcal{D}}^{2p}(X, \mathbb{Z}(p)),
\]
is given by (2.6) $\lambda: H^{p}(X, \overline{K}^{M}_{p}(O_{X}))  \to \mathbb{H}^{2p}(X, \overline{\mathbb{Z}}(p)_{\mathcal{D}}) $, which is the composition:
\[
H^{p}(X, \overline{K}^{M}_{p}(O_{X})) \to H^{p}(\Omega_{X/ \mathbb{C}}^{p-1}) \hookrightarrow \mathbb{H}^{2p}(X, \overline{\mathbb{Z}}(p)_{\mathcal{D}}) .
\]
\end{lemma}

\begin{remark}
 When we consider the infinitesimal deformation of the Deligne complex $\mathbb{Z}(p)_{\mathcal{D}}$, $\mathbb{Z}(p)$ is fixed so that it does not appear in the tangent complex $\overline{\mathbb{Z}}(p)_{\mathcal{D}}$. This explains why the construction of $\lambda$ (2.6) is simpler than the known construction \cite{El-Z, EV, Jannsen} of the Deligne cycle class map.
 
\end{remark}

In the following, we consider $H^{p}(X, \overline{K}^{M}_{p}(O_{X})) \to H^{p}(\Omega_{X/ \mathbb{C}}^{p-1})$ as the infinitesimal form of the Deligne cycle class map $r: CH^{p}(X) \to H_{\mathcal{D}}^{2p}(X, \mathbb{Z}(p))$, and describe it explicitly.

It is well known that $\overline{K}^{M}_{p}(O_{X}) \xrightarrow{\cong}  \Omega_{X/ \mathbb{Q}}^{p-1}$ and the isomorphism is given by 
\begin{equation}
\{1+ \varepsilon \dfrac{g_{1}}{f_{1}}, f_{2}, \cdots, f_{p} \}  \to (-1)^{p-1}\dfrac{ g_{1}}{f_{1}} \dfrac{d_{\mathbb{Q}}f_{2}}{f_{2}} \wedge \cdots \wedge  \dfrac{d_{\mathbb{Q}}f_{p}}{f_{p}}.
\end{equation}
Moreover, one has the following commutative diagram,
\[
\begin{CD}
\{1+ \varepsilon \dfrac{g_{1}}{f_{1}}, f_{2}, \cdots, f_{p} \}  @> \beta >(2.4)>   (-1)^{p-1}\dfrac{ g_{1}}{f_{1}} \dfrac{d_{\mathbb{C}}f_{2}}{f_{2}} \wedge \cdots \wedge  \dfrac{d_{\mathbb{C}}f_{p}}{f_{p}} \\
     @V \cong V(2.7)V   @VV=V \\
 (-1)^{p-1} \dfrac{ g_{1}}{f_{1}}  \dfrac{d_{\mathbb{Q}}f_{2}}{f_{2}} \wedge \cdots \wedge  \dfrac{d_{\mathbb{Q}}f_{p}}{f_{p}}   @>d_{\mathbb{Q}} \to d_{\mathbb{C}}>>  (-1)^{p-1} \dfrac{ g_{1}}{f_{1}} \dfrac{d_{\mathbb{C}}f_{2}}{f_{2}} \wedge \cdots \wedge  \dfrac{d_{\mathbb{C}}f_{p}}{f_{p}}.
    \end{CD}
\]

This shows the following diagram is commutative,
\[
\begin{CD}
       \overline{K}^{M}_{p}(O_{X}) @>\beta>(2.4)>  \Omega_{X/ \mathbb{C}}^{p-1}\\
     @V \cong V(2.7)V   @VV=V \\
\Omega_{X/ \mathbb{Q}}^{p-1} @>d_{\mathbb{Q}} \to d_{\mathbb{C}}>>   \Omega_{X/ \mathbb{C}}^{p-1}.
  \end{CD}
\]

Passing to cohomology groups, one has the following commutative diagram,
\[
\begin{CD}
       H^{p}(X, \overline{K}^{M}_{p}(O_{X}))  @>>>  H^{p}(\Omega_{X/ \mathbb{C}}^{p-1})\\
     @V \cong VV   @VV=V \\
    H^{p}(\Omega_{X/ \mathbb{Q}}^{p-1}) @>d_{\mathbb{Q}} \to d_{\mathbb{C}}>>   H^{p}(\Omega_{X/ \mathbb{C}}^{p-1}).
  \end{CD}
\]

To summarize,
\begin{theorem} \label{theorem: inf'l form2}
Let $X$ be a smooth projective variety defined over $\mathbb{C}$, for each positive integer $p$, the infinitesimal form of the Deligne cycle class map
\[
r: CH^{p}(X) \to H_{D}^{2p}(X, \mathbb{Z}(p)),
\]
is given by 
\[
\delta r: H^{p}(\Omega_{X/ \mathbb{Q}}^{p-1}) \to H^{p}(\Omega_{X/ \mathbb{C}}^{p-1}),
\]
where $\delta r$ is induced by the natural map $\Omega_{X/ \mathbb{Q}}^{p-1} \to \Omega_{X/ \mathbb{C}}^{p-1}$.
\end{theorem}

Let $X$ be a smooth projective variety over $k$, where $k$ is a field of characteristic $0$. For each positive integer $p$, one still has the Soul\'e's variant of Bloch-Quillen identification 
\begin{equation}
 CH^{p}(X)= H^{p}(X, K^{M}_{p}(O_{X})) \  \mathrm{ modulo \ torsion},
\end{equation}
where $K^{M}_{p}(O_{X})$ is the 
Milnor K-theory sheaf associated to the presheaf
\[
  U \to K^{M}_{p}(O_{X}(U)).
\]
The formal tangent space to $CH^{p}(X)$ is identified with $H^{P}(\Omega_{X/ \mathbb{Q}}^{p-1})$.

\begin{corollary}
Let $X$ be a smooth projective variety defined over $k$, where $k$ is a field of characteristic $0$, for each positive integer $p$, the infinitesimal form of the composition
\[
CH^{p}(X)_{\mathbb{Q}} \rightarrow CH^{p}(X_{\mathbb{C}})_{\mathbb{Q}} \xrightarrow{r}  H_{\mathcal{D}}^{2p}(X_{\mathbb{C}}, \mathbb{Z}(p))_{\mathbb{Q}},
\]
has the form
\[
H^{p}(\Omega_{X/ \mathbb{Q}}^{p-1}) \to H^{p}(\Omega_{X_{\mathbb{C}}/ \mathbb{Q}}^{p-1}) \xrightarrow{\delta r} H^{p}(\Omega_{X_{\mathbb{C}}/ \mathbb{C}}^{p-1}),
\]
which is
\begin{equation}
H^{p}(\Omega_{X/ \mathbb{Q}}^{p-1}) \to H^{p}(\Omega_{X_{\mathbb{C}}/ \mathbb{C}}^{p-1}).
\end{equation}
\end{corollary}

If $k$ is a number field, $\Omega_{X/ \mathbb{Q}}=\Omega_{X/k}$ and 
$H^{p}(\Omega_{X/ \mathbb{Q}}^{p-1}) = H^{p}(\Omega_{X/k}^{p-1})$.
By base change, $H^{p}(\Omega_{X_{\mathbb{C}}/ \mathbb{C}}^{p-1}) \cong H^{p}(\Omega_{X/ k}^{p-1}) \otimes_{k} \mathbb{C}$.
The map (2.9) can be rewritten as 
\[
H^{p}(\Omega_{X/ k}^{p-1}) \to H^{p}(\Omega_{X / k}^{p-1}) \otimes_{k} \mathbb{C},
\]
which is obviously injective. In summary,
\begin{theorem} \label{theorem: main theorem}
The infinitesimal form of Conjecture \ref{conjecture: main} is true. To be precise,
let $X$ be a smooth projective variety defined over a number field $k$, for each positive integer $p$, the infinitesimal form of the composition
\[
CH^{p}(X)_{\mathbb{Q}} \rightarrow CH^{p}(X_{\mathbb{C}})_{\mathbb{Q}} \to  H_{\mathcal{D}}^{2p}(X_{\mathbb{C}}, \mathbb{Z}(p))_{\mathbb{Q}}
\]
 in Conjecture \ref{conjecture: main} is of the form, 
\[
H^{p}(\Omega_{X/ k}^{p-1}) \to H^{p}(\Omega_{X/ k}^{p-1}) \otimes_{k} \mathbb{C},
\]
which is injective.

\end{theorem}

Conjecture \ref{conjecture: main} is part of the Bloch-Beilinson conjecture.
Let $X$ be a smooth projective variety over $\mathbb{C}$, for each positive integer $p$,
Bloch-Beilinson conjecture predicts that there is a filtration which has the form
\begin{align*}
CH^{p}(X)_{\mathbb{Q}} & =F^{0}CH^{p}(X)_{\mathbb{Q}} \supset F^{1}CH^{p}(X)_{\mathbb{Q}} \\
& \supset \cdots \supset F^{p}CH^{p}(X)_{\mathbb{Q}} \supset F^{p+1}CH^{p}(X)_{\mathbb{Q}}=0.
\end{align*}
The first two steps are known and $F^{2}CH^{p}(X)_{\mathbb{Q}}$ is the kernel of the Deligne cycle class map 
\[
r: CH^{p}(X)_{\mathbb{Q}} \to  H_{\mathcal{D}}^{2p}(X, \mathbb{Z}(p))_{\mathbb{Q}}.
\]

An alternative way to state Conjecture \ref{conjecture: main} is,
\begin{conjecture} [ cf. Implication 1.2 \cite{GG} page 478]
If $X$ is a smooth projective variety over a number field $k$, then
\[
F^{2}CH^{p}(X)_{\mathbb{Q}}=0,
\]
where $F^{2}CH^{p}(X)_{\mathbb{Q}}$ is the filtration induced from $F^{2}CH^{p}(X_{\mathbb{C}})_{\mathbb{Q}}$ under the natural map $CH^{p}(X) \to CH^{p}(X_{\mathbb{C}})$.
\end{conjecture}

Theorem \ref{theorem: main theorem} suggests that this Conjecture(and Conjecture \ref{conjecture: main}) looks reasonable at the infinitesimal level.

\begin{remark}
The assumption``{\it $k$ is a number field}" in Theorem \ref{theorem: main theorem} is crucial, it guarantees $\Omega_{X/ \mathbb{Q}}=\Omega_{X/k}$.
This suggests that the assumption``{\it $k$ is a number field}" in Conjecture \ref{conjecture: main} can not be loosened. In fact, for the ground field of transcendental degree one, Green-Griffiths-Paranjape \cite{GGP} has found
counterexamples, extending earlier examples by Bloch, Nori and Schoen.

To understand algebraic cycles, the transcendental degree of the ground field does matter.
\end{remark} 

\textbf{Acknowledgements}. This note is inspired by H\'el\`ene Esnault's talk on
\cite{EH} at Tsinghua University (December 2017). The author sincerely thanks Spencer Bloch, H\'el\`ene Esnault, Jerome Hoffman and Jan Stienstra for discussions. Jerome Hoffman has read a preliminary version, his suggestions improve this note a lot. This work is partially supported by the Fundamental Research Funds for the Central Universities.

\end{document}